\documentclass[11pt]{article}
\usepackage{amssymb, amsmath,fullpage}
\numberwithin{equation}{section}
\newtheorem {thm}{Theorem}[section]
\newtheorem {lem}[thm]{Lemma}
\newtheorem {cor}[thm]{Corollary}
\newtheorem {prop}[thm]{Proposition}

\newcommand{\la}{\langle}
\newcommand{\lla}{\langle\!\langle}
\newcommand{\ra}{\rangle}
\newcommand{\rra}{\rangle\!\rangle}
\newcommand{\proof}{\par\noindent{\em Proof.\ }}

\newcommand{\qed}{\hfill $\square$\par\smallskip}
\renewcommand{\vec}[1]{\boldsymbol{#1}}
\title{\bf\boldmath 
A generalization of the Sears--Slater transformation \\
and elliptic Lagrange interpolation of type $BC_n$ 
}
\author{
{\sc Masahiko Ito}\thanks{
Corresponding Author at: School of Science and Technology for Future Life,
Tokyo Denki University,
Tokyo 120-8551, Japan, E-mail: 
{\tt mito@cck.dendai.ac.jp}}\ \ 
and {\sc Masatoshi Noumi}\thanks{
 Department of Mathematics, Kobe University, 
 Rokko, Kobe 657-8501, Japan, E-mail: 
{\tt noumi@math.kobe-u.ac.jp}
 }
}
\date{%June 23, 2015
}

\begin{document}
\maketitle
\begin{abstract}
The connection formula for the Jackson integral of type $BC_n$ is obtained in the form of 
a Sears--Slater type expansion of a bilateral multiple basic hypergeometric series as 
a linear combination of several specific bilateral multiple series. The coefficients of this expansion are expressed by certain elliptic Lagrange interpolation functions.  
Analyzing basic properties of  
the elliptic Lagrange interpolation functions, 
an explicit determinant formula is provided 
for a fundamental solution matrix of the associated system of $q$-difference equations.

\end{abstract}



\section{Introduction}
Throughout this paper 
we fix the base $q\in\mathbb{C}$ with $0<|q|<1$, 
and use the notation of $q$-shifted factorials 
\begin{equation*}
\begin{split}
&(u)_\infty=\prod^\infty_{l=0}(1-uq^l),\quad
(a_1,\dots,a_r)_\infty=(a_1)_\infty\cdots(a_r)_\infty,
\\
&(u)_\nu = (u)_\infty/(uq^\nu)_\infty,\quad
(a_1,\dots,a_r)_\nu=(a_1)_\nu\cdots(a_r)_\nu\quad 
(\nu\in\mathbb{Z}). 
\end{split}
\end{equation*}
For generic complex parameters $a_1,\ldots,a_r$ and 
$b_1,\ldots,b_r$, 
the bilateral basic hypergeometric series $_r\psi_r$ is defined by 
\begin{equation*}
 \label{eq:6.1}
  {}_r\psi_r
  \left[\!\!
\begin{array}{c}
      a_1,\dots,a_r \\
      b_1,\dots,b_r
\end{array}\!; q,x
  \right]
= \sum_{\nu=-\infty}^\infty\frac{(a_1,\dots,a_r)_\nu}{(b_1,\dots,b_r)_\nu}x^\nu\qquad
(|b_1\cdots b_r/a_1\cdots a_r|<|x|<1). 
\end{equation*}
There is a celebrated summation formula for 
very well-poised balanced ${}_6\psi_6$ series
\begin{equation}
\label{eq:66}
{}_6\psi_6
\bigg[\!\!
\begin{array}{c}
      q\sqrt{a},-q\sqrt{a},\, b,\, c,\, d,\, e\, \\
      \sqrt{a}, \, -\sqrt{a},\!\frac{aq}{b},\! \frac{aq}{c},\! \frac{aq}{d},\! \frac{aq}{e}
\end{array}\!; q,{\displaystyle\frac{a^2q}{bcde}}
\bigg]=
\frac{(aq,\frac{aq}{bc},\frac{aq}{bd},\frac{aq}{be},\frac{aq}{cd},\frac{aq}{ce},\frac{aq}{de},q,\frac{q}{a})_\infty}
{(\frac{aq}{b},\frac{aq}{c},\frac{aq}{d},\frac{aq}{e},\frac{q}{b},\frac{q}{c},\frac{q}{d},\frac{q}{e},\frac{a^2q}{bcde}
)_\infty}
\end{equation}
for $|a^2q/bcde|<1$, 
called {\it Bailey's sum}
(see \cite[(5.3.1), p.\,140]{GR2004}).
On the other hand, 
there are several transformation 
formulas for well-poised or very well-poised
${}_{2r}\psi _{2r}$ series due to  
Sears \cite{Se1,Se2} and Slater \cite{Sl1,Sl2}.
A typical example is 
\begin{eqnarray}
\label{eq:2r2r}
\lefteqn{
{}_{2r} \psi _{2r}\bigg[
\begin{array}{ccccc}
		\!\!q\sqrt{a},&\!\!\!\! -q\sqrt{a},&\!\!\!\! b_{3},&\!\!\!\!\ldots,&\!\!\!\! b_{2r}\!\!\\
		\!\!\sqrt{a},&\!\!\!\! -\sqrt{a},&\!\!\!\! \frac{aq}{b_{3}},&\!\!\!\!\ldots,&\!\!\!\!\frac{aq}{b_{2r}}\!\!
		\end{array}; q,
		\frac{a^{r-1}q^{r-2}}{b_{3}\ldots b_{2r}}\bigg]}
\nonumber\\[5pt]
&=&\frac{
	(a_{4},\ldots,a_{r},\frac{q}{a_{4}},\ldots,\frac{q}{a_{r}},
	\frac{a_{4}}{a},\ldots,\frac{a_{r}}{a},\frac{aq}{a_{4}},\ldots,\frac{aq}{a_{r}},
	\frac{a_{3}q}{b_{3}},\ldots,\frac{a_{3}q}{b_{2r}},\frac{aq}{a_{3}b_{3}},\ldots,\frac{aq}{a_{3}b_{2r}},aq,\frac{q}{a})_{\infty}
	}
	{(\frac{q}{b_{3}},\ldots,\frac{q}{b_{2r}},\frac{aq}{b_{3}},\ldots,\frac{aq}{b_{2r}},
	\frac{a_{4}}{a_{3}},\ldots,\frac{a_{r}}{a_{3}},\frac{a_{3}q}{a_{4}},\ldots,\frac{a_{3}q}{a_{r}}, 
	\frac{a_{3}a_{4}}{a},\ldots,\frac{a_{3}a_{r}}{a},\frac{aq}{a_{3}a_{4}},\ldots,\frac{aq}{a_{3}a_{r}}
	,\frac{a_{3}^2q}{a},\frac{aq}{a_{3}^2})_{\infty}
	}
\nonumber\\[5pt]
&&\quad\times
\,{}_{2r}\psi _{2r}\bigg[
\begin{array}{ccccc}
		\!\!\frac{qa_{3}}{\sqrt{a}},&\!\!\!\!-\frac{qa_{3}}{\sqrt{a}},&\!\!\!\!\frac{a_{3}b_{3}}{a},&\!\!\!\!\ldots,&\!\!\!\!\frac{a_{3}b_{2r}}{a}\!\!\\
		\!\!\frac{a_{3}}{\sqrt{a}},&\!\!\!\!-\frac{a_{3}}{\sqrt{a}},&\!\!\!\!\frac{a_{3}q}{b_{3}},&\!\!\!\!\ldots,&\!\!\!\!\frac{a_{3}q}{b_{2r}}\!\!
		\end{array};q,
		\frac{a^{r-1}q^{r-2}}{b_{3}\ldots b_{2r}}\bigg]
\nonumber\\[5pt]
&&+\ {\rm idem} (a_{3};a_{4},\ldots,a_{r})%\nonumber
\end{eqnarray}
for $|a^{r-1}q^{r-2}/b_3\cdots b_{2r}|<1$,
which is called {\it Slater's transformation formula 
for a very well-poised balanced ${}_{2r}\psi_{2r}$ series}
(see \cite[(5.5.2), p.\,143]{GR2004}). 
Here the symbol ``idem$(a_{3};a_{4},\ldots,a_{r})$" 
stands for 
the sum of the $r-3$ expressions obtained from the preceding 
one 
by interchanging $a_{3}$ with each $a_{k}$ ($k=4,\ldots,r$).
\\

The $BC_n$ Jackson integrals, 
which we are going to discuss below, are 
a multiple sum generalization of 
the very well-poised 
${}_{2r}\psi_{2r}$ series.  
They provide with a natural framework 
for summation/transformation formulas 
for basic hypergeometric series, 
from the viewpoint of the Weyl group symmetry and the $q$-difference equations.  
The holonomic system of $q$-difference equations satisfied by a $BC_n$ Jackson integral has been investigated 
in \cite{AI2009}, from which we recall some terminology.

For a function $\varphi=\varphi(z)$ of 
$z=(z_1,\ldots,z_n)\in(\mathbb{C}^*)^n$, 
we denote by 
\begin{equation*}
\la\varphi, z\ra=
\int_0^{z\infty} \varphi(w)\Phi(w)\Delta(w)
\frac{d_qw_1}{w_1}\wedge\cdots\wedge
\frac{d_qw_n}{w_n}
=(1-q)^n\sum_{\nu\in \mathbb{Z}^n}\varphi(z q^\nu)\Phi(z q^\nu)\Delta(z q^\nu),
\end{equation*}
the Jackson integral associated with the multiplicative lattice 
$z q^\nu=(z_1q^{\nu_1},\ldots,z_nq^{\nu_n})\in(\mathbb{C}^*)^n$ 
($\nu=(\nu_1,\ldots,\nu_n)\in \mathbb{Z}^n$). 
In this definition, we specify the integrand by the weight function
$$
\Phi(z)=\prod_{i=1}^n\prod_{m=1}^{2s+2}z_i^{\frac{1}{2}-\alpha_m}\frac{(qa_m^{-1}z_i)_\infty}{(a_mz_i)_\infty}
\prod_{1\le j<k\le n}z_j^{1-2\tau}\frac{(qt^{-1}z_j/z_k)_\infty(qt^{-1}z_jz_k)_\infty}{(tz_j/z_k)_\infty(tz_jz_k)_\infty},
$$
where $q^{\alpha_m}=a_m$ and $q^{\tau}=t$, and the Weyl denominator of type $C_n$
$$
\Delta(z)=\prod_{i=1}^n\frac{1-z_i^2}{z_i}
\prod_{1\le j<k\le n}\frac{(1-z_j/z_k)(1-z_jz_k)}{z_j}.
$$
We call the sum $\la\varphi, z\ra$ the {\it Jackson integral of type $BC_n$} if it converges. 
We denote by $W_n=\{\pm1\}^n\rtimes \mathfrak{S}_n$ the Weyl group of type $C_n$ (hyperoctahedral group of degree $n$); 
this group acts on the field of meromorphic functions on 
$({\Bbb C}^*)^n$ 
through the permutations and the inversions of variables $z_1,\ldots,z_n$.  
Setting 
$$
\lla\varphi, z\rra=\frac{\la\varphi, z\ra}{\Theta(z)},\quad
\Theta(z)=\prod_{i=1}^n
\frac{z_i^s\theta(z_i^2)}{\prod_{m=1}^{2s+2}z_i^{\alpha_m}\theta(a_mz_i)}
\prod_{1\le j<k\le n}
\frac{\theta(z_j/z_k)\theta(z_jz_k)}{z_j^{2\tau}\theta(tz_j/z_k)\theta(tz_jz_k)},
$$
where $\theta(u)=(u)_\infty (qu^{-1})_\infty$, 
we call $\lla\varphi, z\rra$ the {\it regularized Jackson integral of type $BC_n$}. 
We remark that, 
for any $W_n$-invariant 
holomorphic function $\varphi(z)$ on $({\Bbb C}^*)^n$, 
the regularization $\lla\varphi, z\rra$ 
is also holomorphic and $W_n$-invariant 
as a function of $z\in(\mathbb{C}^*)^n$ 
(see \cite[Definition 3.8]{AI2009}).
This regularized Jackson integral 
$\lla\varphi, z\rra$ is the main object of this paper. 

In what follows, we set 
$$
B=B_{s,n}=\{\lambda=(\lambda_1,\lambda_2,\ldots, \lambda_n)\in \mathbb{Z}^n\,;\,s-1\ge\lambda_1\ge\lambda_2\ge \cdots\ge \lambda_n\ge 0\},\quad 
$$
so that $|B_{s,n}|={s+n-1\choose n}$. 
We also use the symbol $\preceq$ for the lexicographic order
of $B_{s,n}$.  
Namely,  
for $\lambda,\mu\in B_{s,n}$, we denote $\lambda\prec\mu$ if 
there exists $k\in\{1,2,\ldots,n\}$ such that 
$
\lambda_1=\mu_1,\lambda_2=\mu_2,\ldots,\lambda_{k-1}=\mu_{k-1}
$
and
$\lambda_k<\mu_k$.
For each $\lambda\in B_{s,n}$, 
we denote by $\chi_\lambda(z)$ 
the {\em symplectic Schur function} 
$$
\chi_\lambda(z)=
\frac{\det(z_i^{\lambda_j+n-j+1}-z_i^{-\lambda_j-(n-j+1)})_{1\le i,j\le n}}
{\det(z_i^{n-j+1}-z_i^{-(n-j+1)})_{1\le i,j\le n}}
=\frac{\det(z_i^{\lambda_j+n-j+1}-z_i^{-\lambda_j-(n-j+1)})_{1\le i,j\le n}}
{\Delta(z)},
$$
which is a $W_n$-invariant Laurent polynomial. 
With the basis $\{\chi_\lambda(z)\,;\,\lambda\in B_{s,n}\}$, 
one can construct a holonomic system of 
$q$-difference equations of rank ${s+n-1\choose n}$ 
for the Jackson integral of type $BC_n$. 
\begin{prop}[\boldmath $q$-Difference system \cite{AI2005,AI2008}]
Let $\vec{v}(z)$ 
be row vector of functions in $a_1,\ldots,a_{2s+2}$ and $t$ 
specified by 
$$\vec{v}(z)=\big(\lla\chi_\lambda,z\rra\big)_{\lambda\in B},$$ 
where the indices $\lambda\in B_{s,n}$ 
are arranged in the increasing order by $\preceq$.
If the parameters $a_1,a_2,\ldots,a_{2s+2}$ and $t$ are generic, 
there exist invertible ${s+n-1\choose n}\times {s+n-1\choose n}$ matrices 
$Y_{a_i} (i = 1,2,\ldots,2s + 2)$ and $Y_t$ whose entries are rational functions of $a_1,a_2,\ldots,a_{2s+2}$ and $t$, 
not depending on $z$, such that
\begin{equation}
\label{eq:system}
T_{a_i}\vec{v}(z)=\vec{v}(z)\, Y_{a_i},\qquad T_{t}\,\vec{v}(z)=\vec{v}(z)\, Y_{t},
\end{equation}
where $T_{u}$ stands for the $q$-shift operator 
in $u$.
\end{prop}

To construct independent solutions of the above system, 
we define ${s+n-1\choose n}$ specific points in $(\mathbb{C}^*)^n$ as follows. 
Setting
$$
Z=Z_{s,n}=\{\mu=(\mu_1,\mu_2,\ldots,\mu_s)\in \mathbb{N}^s\,;\,\mu_1+\mu_2+\cdots+\mu_s=n\}
$$
so that $|Z_{s,n}|={s+n-1\choose n}$, we denote 
the lexicographic order of $Z_{s,n}$ by $\preceq$. 
For an arbitrary $x=(x_1,x_2,\ldots,x_s)\in (\mathbb{C}^*)^s$, 
we consider the points $x_\mu$ $(\mu \in Z_{s,n})$ 
in $(\mathbb{C}^*)^n$ specified by 
\begin{equation}
\label{eq:x_mu}
x_\mu=(
\underbrace{x_1,x_1t,\ldots,x_1t^{\mu_1-1}\phantom{\Big|}\!\!}_{\mu_1}
,\underbrace{x_2,x_2t,\ldots,x_2t^{\mu_2-1}\phantom{\Big|}\!\!}_{\mu_2}
,\ldots
,\underbrace{x_s,x_st,\ldots,x_st^{\mu_s-1}\phantom{\Big|}\!\!}_{\mu_s})
\in (\mathbb{C}^*)^n.
\end{equation}
In order to confirm that 
the solutions $\vec{v}(x_{\mu})$ ($\mu\in Z_{s,n}$) of the system  (\ref{eq:system}) are linearly independent, 
we need to verify that the determinant of the matrix 
$\big(\lla\chi_\lambda,x_{\mu}\rra\big)_{\lambda\in B,\mu\in Z}$ 
(``Wronskian'' of the system (\ref{eq:system})) is nonzero 
under the genericity condition for parameters. 
One of our main results is the following.
\begin{thm}[Determinant formula]\label{thm:Wronski}
The determinant of the matrix 
$\big(\lla\chi_\lambda,x_{\mu}\rra\big)_{\lambda\in B,\mu\in Z}$ 
is represented explicitly as 
\begin{eqnarray}
\label{eq:Wronski}
&&
\hspace{-16pt}
\det\Big(\lla\chi_\lambda,x_{\mu}\rra\Big)_{\!\!\lambda\in B\atop \!\!\mu\in Z}
=
\prod_{k=1}^n\bigg[\bigg(\!
(1-q)\frac{(q)_\infty(qt^{-(n-k+1)})_\infty}{(qt^{-1})_\infty}\!\bigg)^{\!\! s}\,
\frac{\prod_{1\le i<j\le 2s+2}(qt^{-(n-k)}a_i^{-1}a_j^{-1})_\infty}{(qt^{-(n+k-2)}a_1^{-1}a_2^{-1}\cdots a_{2s+2}^{-1})_\infty}
\bigg]^{{s+k-2\choose k-1}}
\nonumber\\
&&
\hspace{90pt}
\times\prod_{k=1}^n\bigg[
\prod_{r=0}^{n-k}\prod_{1\le i<j\le s}
\frac{\theta(t^{2r-(n-k)}x_ix_j^{-1})\theta(t^{n-k}x_ix_j)}{t^rx_i}
\bigg]^{{s+k-3\choose k-1}},
\end{eqnarray}
where 
the rows $\lambda\in B$ and the columns $\mu\in Z$ of the matrix are arranged by $\preceq$, respectively.
\end{thm}

Notice that our explicit formula \eqref{eq:Wronski}
of the determinant 
splits into two parts. 
The first part is independent of the choice of cycles 
of the integral 
(i.e., independent of $x$), while the second is 
a function of $x$ only.  
From this fact we immediately see that the determinant 
does not vanish if $x$ is generic.

We remark that in the case $s=1$
the matrix size of $(\lla\chi_\lambda,x_\mu\rra)_{\lambda,\mu}$ reduces to 1 
and (\ref{eq:Wronski}) in Theorem \ref{thm:Wronski} becomes 
the following formula first proved by van Diejen \cite{vD1997}:
\begin{equation}
\label{eq:vD-sum}
\lla 1,z\rra=
\prod_{k=1}^n
(1-q)\frac{(q)_\infty(qt^{-k})_\infty}{(qt^{-1})_\infty}
\frac{\prod_{1\le i<j\le 4}(qt^{-(n-k)}a_i^{-1}a_j^{-1})_\infty}{(qt^{-(n+k-2)}a_1^{-1}a_2^{-1}a_3^{-1} a_4^{-1})_\infty},
\end{equation}
which is equivalent to the $q$-Macdonald--Morris identity of type $(C_n^\vee,C_n)$ studied by Gustafson \cite{Gu1990}. 
(See \cite{Ito2006} for the derivation of (\ref{eq:vD-sum}) along the context of this paper. 
See also \cite{Ito2006-2,Kad} for the other derivations.) 
In this case the last factor including theta functions in (\ref{eq:Wronski}) disappears. 
Since (\ref{eq:vD-sum}) coincides with (\ref{eq:66}) if $n=1$, 
we can regard (\ref{eq:Wronski}) as a further extension of Bailey's ${}_6\psi_6$ summation theorem. 
\\

Next we consider the connection problem among the independent cycles.  
From a property of functions written by integral representation, 
the holonomic system of $q$-difference equations satisfied by $\lla \varphi,z\rra$ does not depend on 
the choice of cycles, i.e. the choice of points 
$z\in (\mathbb{C}^*)^n$. 
It also turns out that 
$\lla \varphi, z\rra$ for any $z\in (\mathbb{C}^*)^n$ is expressed as a linear combination of the special solutions $\lla\varphi, x_\mu\rra$, $\mu\in Z_{s,n}$, of the system, i.e.,
\begin{equation}
\label{eq:connection1}
\lla \varphi, z\rra=\sum_{\mu\in Z}c_\mu \lla\varphi, x_\mu\rra, 
\end{equation}
where $c_\mu$ are some connection coefficients. 
We introduce some terminology 
before we state the explicit form of the connection coefficients 
$c_\mu$. 

Let ${\cal O}(({\Bbb C}^*)^n)$ be the $\mathbb{C}$-vector 
space of holomorphic functions on $({\Bbb C}^*)^n$. 
We consider the $\mathbb{C}$-linear subspace 
$H_{m,n}\subset  {\cal O}(({\Bbb C}^*)^n$ 
consisting of all $W_n$-invariant holomorphic functions $f(z)$ such that $T_{z_i}f(z)=f(z)/(qz_i^2)^m\ (i=1,\ldots,n)$, where $T_{z_i}$ stands for the $q$-shift operator 
in $z_i$:
\begin{equation}
\label{eq:def-Hmn}
H_{m,n}=\{f(z)\in\mathcal{O}((\mathbb{C}^*)^n)^{W_n}\,;\,
T_{z_i}f(z)=(qz_i^2)^{-m}f(z)\ \ (i=1,\ldots,n)\}.
\end{equation}
The dimension of $H_{m,n}$ as a $\mathbb{C}$-vector space is known to be $m+n\choose n$. 
(See \cite[Lemma 3.2]{Ito2008}.)
\begin{thm}
\label{thm:E}
For generic $x\in (\mathbb{C}^*)^s$ 
there exists a unique basis $\{f_\lambda(z)\,;\,\lambda\in Z_{s,n}\}$ of $H_{s-1,n}$ such that 
\begin{equation}
\label{eq:E-delta}
f_\lambda(x_\mu)=\delta_{\lambda\mu}\qquad
(\lambda,\mu\in Z_{s,n}),
\end{equation}
where $\delta_{\lambda\mu}$ 
is the Kronecker delta. We denote the function $f_\lambda(z)$ by the symbol $E_\lambda(x;z)$. 
\end{thm}

We call $E_\lambda(x;z)$ the {\it elliptic Lagrange interpolation 
functions of type $BC_n$}. 
An explicit construction of these functions $E_\lambda(x;z)$ 
will be given 
In Section \ref{Section:2}.  
In particular we will prove 
\begin{thm}
\label{thm:E-explicit}
The elliptic Lagrange interpolation functions $E_\lambda(x;z)$ 
of type $BC_n$ are represented explicitly as 
\begin{equation}
\label{eq:E-explicit}
E_\lambda(x;z)=\sum_{K_1\sqcup\cdots\sqcup K_s\atop =\{1,2,\ldots,n\}}
\prod_{i=1}^s\prod_{k\in K_i}
\prod_{1\le j\le s\atop j\ne i}
\frac{\theta(x_jt^{\lambda_j^{(k-1)}}z_{k})\theta(x_jt^{\lambda_j^{(k-1)}}z_{k}^{-1})}
{\theta(x_jt^{\lambda_j^{(k-1)}}x_it^{\lambda_i^{(k-1)}})\theta(x_jt^{\lambda_j^{(k-1)}}x_i^{-1}t^{-\lambda_i^{(k-1)}})},
\end{equation}
where $\lambda_i^{(k)}=|K_i\cap\{1,2,\ldots,k\}|$, 
and the summation is taken over all 
partitions $K_1\sqcup\cdots\sqcup K_s=\{1,2,\ldots,n\}$ 
such that 
$|K_i|=\lambda_i$ $(i=1,2,\ldots,s)$. 
\end{thm}

Here we mention some special cases of Theorem \ref{thm:E-explicit}. \\
\noindent
{\bf Example 1.} The case $n=1$. We have $Z_{s,1}=\{\epsilon_1,\epsilon_2,\ldots,\epsilon_s\}$, 
where $\epsilon_i=(0,\ldots,0,\overset{\text{\tiny$i$}}{\stackrel{\text{\tiny$\smile$}}{1}},0,\ldots,0)$,
and $x_{\epsilon_i}=x_i$, i.e., we obtain 
\begin{equation}
\label{eq:E-explicit,n=1}
E_{\epsilon_i}(x;z)=\prod_{1\le j\le s\atop j\ne i}
\frac{\theta(x_jz)\theta(x_j/z)}{\theta(x_ix_j)\theta(x_j/x_i)}. 
\end{equation}
%\noindent
%{\bf Example 2.} The case $s=1$. We have $Z_{1,n}=\{(n)\}$, which has only one element, 
%and have $x_{(n)}=(x,xt,\ldots,xt^{n-1})$ for $x\in \mathbb{C}^*$. 
%Then $E_{(n)}(x;z)=1$, so that $\lla \varphi, z\rra=\lla\varphi, x_{(n)}\rra$ for any $z\in \mathbb{C}^*$.
%This means that $\lla \varphi, z\rra$ is a constant independent of $z$. In particular, 
%$\lla 1, z\rra=\lla 1, x_{(n)}\rra$ is evaluated as (\ref{eq:vD-sum}).\\
%
\noindent
{\bf Example 2.} The case $s=2$. We have $Z_{2,n}=\{(r,n-r)\,;\, r=0,1,\ldots,n\}$. Then 
\begin{eqnarray*}
E_{(r,n-r)}(x_1,x_2;z)&=&
\sum_{1\le i_1<\cdots<i_r\le n \atop 1\le j_1<\cdots<j_{n-r}\le n}
\prod_{k=1}^r
\frac{\theta(x_2t^{i_k-k}z_{i_k})\theta(x_2t^{i_k-k}z_{i_k}^{-1})}
{\theta(x_2t^{i_k-k}x_1t^{k-1})\theta(x_2t^{i_k-k}x_1^{-1}t^{-(k-1)})}\nonumber\\[-10pt]
&&\hspace{70pt}\times 
\prod_{l=1}^{n-r}
\frac{\theta(x_1t^{j_l-l}z_{j_l})\theta(x_1t^{j_l-l}z_{j_l}^{-1})}
{\theta(x_1t^{j_l-l}x_2t^{l-1})\theta(x_1t^{j_l-l}x_2^{-1}t^{-(l-1)})},
\end{eqnarray*}
where the summation is taken over all pairs of sequences 
$1\le i_1<\cdots<i_r\le n$ and $1\le j_1<\cdots<j_{n-r}\le n$ such that 
$\{i_1,\ldots,i_r\}\cup\{j_1,\ldots,j_{n-r}\}=\{1,2,\ldots,n\}$.
In our previous works \cite{IN1,IN2}, 
these functions $E_{(r,n-r)}(x_1,x_2;z)$ ($0\le r\le n$) 
are called the {\em $BC_n$ fundamental invariants}, 
and are effectively used in 
an alternative method of evaluation of 
the $BC_n$ elliptic Selberg integral proposed 
by van Diejen and Spiridonov \cite{vDS2001}. 
\\

We now return to the connection problem.  
Another main result of this paper is that 
the connection coefficient $c_\mu$ in (\ref{eq:connection1}) exactly coincides with our 
elliptic Lagrange interpolation function $E_\mu(x;z)$ 
for each $\mu\in Z_{s,n}$. 
\begin{thm}[Connection formula] 
\label{thm:connection2}
Suppose that $\varphi(z)\in {\cal O}(({\Bbb C}^*)^n)$ is 
$W_n$-invariant. Then 
\begin{equation}
\label{eq:connection2}
\lla \varphi, z\rra=\sum_{\mu\in Z_{s,n}}\lla\varphi, x_\mu\rra E_\mu(x;z),
\end{equation}
where the connection coefficients $E_\mu(x;z)$ are explicitly written as {\rm (\ref{eq:E-explicit})}.
\end{thm}

We call this connection formula the {\it generalized Sears--Slater transformation}.  In fact, the connection 
formula (\ref{eq:connection2}) of the case $n=1$ is given by
$$
\lla \varphi, z\rra=\sum_{i=1}^s\lla\varphi, x_i\rra 
\prod_{1\le j\le s\atop j\ne i}
\frac{\theta(x_jz)\theta(x_j/z)}{\theta(x_ix_j)\theta(x_j/x_i)},
$$
which exactly coincides with the transformation formula (\ref{eq:2r2r}) if $\varphi(z)\equiv1$ and $r=s+2$. 
See \cite{IS} for details about the correspondence between them. 

\par\medskip
This paper is organized as follows. 
We first provide in Section \ref{Section:2}
a proof of Theorem \ref{thm:E}
based on 
an explicit construction 
of the elliptic Lagrange interpolation functions 
by means of 
a kernel function as in \cite{KNS}. 
Section \ref{Section:3} is devoted to proving 
the explicit formula of Theorem \ref{thm:E-explicit}
for our interpolation functions.  
In Section \ref{Section:4} we investigate the transition coefficients between two sets of the interpolation functions. 
In particular we propose an explicit formula for the determinant 
of the transition matrix. 
Using properties of the elliptic Lagrange interpolation 
functions, we complete 
in Section \ref{Section:5} 
the proof of the connection formula 
of Theorem \ref{thm:connection2}.  
The determinant formula 
of Theorem \ref{thm:Wronski} is also obtained 
as a corollary of Theorem \ref{thm:connection2}. 

Lastly we remark that  
another type of $BC_n$ 
Jackson integral is studied in \cite{Ito2008}.  
In that case, the determinant formula 
\cite[Theorem 1.7]{Ito2008} corresponding 
to Theorem \ref{thm:Wronski} is 
regarded as a generalization of Gustafson's $C_n$ sum.  
The corresponding connection formula \cite[Theorem 1.1]{Ito2008}
is much simpler than 
Theorem \ref{thm:connection2} of this paper.  

%%%%%%%%%%%%%%%%%%%%%%%%%%%%%%%%%%
%%%%%%%%%%%%%%%%%%%%%%%%%%%%%%%%%%
\section{\boldmath 
Construction of the $BC_n$ interpolation functions}
%%%%%%%%%%%%%%%%%%%%%%%%%%%%%%%%%%
%%%%%%%%%%%%%%%%%%%%%%%%%%%%%%%%%%
\label{Section:2}
%%%%%%%%%%%%%%%%%%%%%%%%%%%%%%%%%%
In this section we give a proof of Theorem \ref{thm:E}.\\

Throughout this paper we use the symbol 
$$
e(a;b)=a^{-1}\theta(ab)\theta(ab^{-1})\qquad
(a,b\in\mathbb{C}^\ast), 
$$ 
where $\theta(u)=(u)_\infty(q/u)_\infty$. 
Since $\theta(a)=\theta(qa^{-1})$ and $\theta(qa)=-a^{-1}\theta(a)$, 
this symbol satisfies 
$$
e(a^{-1};b)=e(a;b),\quad e(a;b)=-e(b;a),
\quad e(a;a)=0
\quad\mbox{and}\quad e(qa;b)=
(qa^2)^{-1}e(a;b). 
$$
Fixing a generic parameter $t\in\mathbb{C}^\ast$, we also introduce 
the notation 
of {\em $t$-shifted factorials} 
$$
e(a;b)_r=e(a;b)e(at;b)\cdots e(at^{r-1};b)
\quad(r=0,1,2,\ldots)
$$
associated with the symbol $e(a;b)$.  
\par\medskip
For two sets of variables 
$z\in (\mathbb{C}^*)^{n}$ and 
$y\in (\mathbb{C}^*)^{s-1}$, 
we consider the {\em dual Cauchy kernel} 
$$\Psi(z;y)=\prod_{i=1}^n\prod_{j=1}^{s-1}e(z_i;y_j).$$ 
Note that $\Psi(z,y)$ is a holomorphic function 
on $(\mathbb{C}^{\ast})^{n}\times (\mathbb{C}^\ast)^{s-1}$, 
and satisfies 
$$\Psi(z;y)\in H_{s-1,n}^z\quad\mbox{and}\quad 
\Psi(z;y)\in H_{n,s-1}^y,$$
with superscripts indicating the variables. 
For each multi-index $\mu\in Z_{s,n}$, 
we define a function $F_{\mu}(x;y)$ of $y\in(\mathbb{C}^\ast)^{s-1}$ with parameters $x\in(\mathbb{C}^\ast)^s$ by
\begin{equation}
\label{eq:def-F}
F_\mu(x;y)
=\Psi(x_\mu;y)=\prod_{i=1}^{s}\prod_{j=1}^{s-1}e(x_i;y_j)_{\mu_i}
\quad(x\in(\mathbb{C}^\ast)^s, \ y\in
(\mathbb{C}^\ast)^{s-1}), 
\end{equation}
where $x_\mu\in (\mathbb{C}^*)^n$ is 
specified by (\ref{eq:x_mu}). 
These functions $F_\mu(x;y)$ $(\mu\in Z_{s,n})$ 
are $W_{s-1}$-invariant with respect to $y$, and 
satisfy
$T_{y_i}F_\mu(x;y)=(qy_i^2)^{-n}F_\mu(x;y)$
$(i=1,\ldots,s-1)$, 
namely, 
$F_\mu(x;y)\in H_{n,s-1}^y$. 
On the other hand, 
for each $x\in(\mathbb{C}^\ast)^s$ and 
$\nu=(\nu_1,\nu_2,\ldots,\nu_s)\in Z_{s,n}$, 
we specify the point $\eta_{\nu}(x)$ in $(\mathbb{C}^*)^{s-1}$ by 
$$
\eta_{\nu}(x)=(x_1t^{\nu_1},x_2t^{\nu_2},\ldots,x_{s-1}t^{\nu_{s-1}})\in (\mathbb{C}^*)^{s-1}.
$$
Note that  
the point $\eta_{\nu}(x)$ has no coordinate corresponding to 
the 
$s$th index $\nu_s$ of $\nu\in Z_{s,n}$, while 
$\eta_{\nu}(x)$ is defined injectively from $\nu\in Z_{s,n}$
if  $x$ is generic.  
In fact, if we set $$L_{s,n}=\{(\nu_1,\nu_2,\ldots,\nu_{s-1})\in \mathbb{N}^{s-1}\,;\,
\nu_1+\nu_2+\ldots+\nu_{s-1}\le n\},$$ then 
the map 
$(\nu_1,\nu_2,\ldots,\nu_{s-1},\nu_s)\mapsto  (\nu_1,\nu_2,\ldots,\nu_{s-1})$ 
defines a bijection $Z_{s,n}\to L_{s,n}$. 

\begin{lem}[Triangularity]
\label{lem:triangular-F}
For each $\mu,\nu\in Z_{s,n}$, 
$F_\mu(x;\eta_{\nu}(x))=0$
unless $\mu_i\le \nu_i$ $(i=1,\ldots,s-1)$.  
In particular, $F_\mu(x;\eta_\nu(x))=0$ 
for $\mu\succ \nu$. 
Moreover, if $x\in (\mathbb{C}^*)^s$ is generic, 
then 
$F_\mu(x;\eta_\mu(x))\ne 0$ for all $\mu\in Z_{s,n}$. 
\end{lem}

This lemma implies that the matrix 
$F=\big(F_\mu(x;\eta_\nu(x))\big)_{\mu,\nu\in Z}$ 
is upper triangular, and also invertible if 
$x\in\mathbb{C}^{s}$ is generic. 
\\

\proof If there exists $j\in \{1,2,\ldots,s-1\}$ such that $\nu_j<\mu_j$, then 
$F_\mu(x;\eta_{\nu}(x))=0$. In fact, 
in the expression 
$$
F_\mu(x;\eta_{\nu}(x))
=\prod_{i=1}^{s}\prod_{j=1}^{s-1}e(x_i;x_jt^{\nu_j})_{\mu_i}, 
$$
the function $e(x_j;x_jt^{\nu_j})_{\mu_j}$ 
has the factor 
$\theta(t^{-\nu_j})\theta(t^{-\nu_j+1})\cdots 
\theta(t^{-\nu_j+(\mu_j-1)})=0$ 
if $\nu_j<\mu_j$. 
If $\nu\prec\mu$,  then $\nu_i<\mu_i$ for some $i\in \{1,2,\ldots,s-1\}$ by definition,  and hence 
we obtain 
$F_\mu(x;\eta_{\nu}(x))=0$ if $\nu\prec\mu$.
For $\nu=\mu$, 
\begin{equation*}
\begin{split}
F_\mu(x;\eta_{\mu}(x))
&=
\prod_{i=1}^{s}\prod_{j=1}^{s-1}
\frac{\theta(x_ix_jt^{\mu_j})}{x_i}
\frac{\theta(x_ix_jt^{\mu_j+1})}{x_it}\cdots 
\frac{\theta(x_ix_jt^{\mu_j+(\mu_i-1)})}{x_it^{\mu_i-1}}
\\
&\qquad\qquad
\times 
\theta(x_ix_j^{-1}t^{-\mu_j})\theta(x_ix_j^{-1}t^{-\mu_j+1})\cdots 
\theta(x_ix_j^{-1}t^{-\mu_j+(\mu_i-1)})
\end{split}
\end{equation*}
does not vanish if we impose an appropriate 
genericity condition on $x\in (\mathbb{C}^*)^s$. 
\qed

\begin{lem}
\label{lem:basis-F}
The set $\{F_\mu(x;y)\,;\, \mu\in Z_{s,n}\}$ is a basis of the $\mathbb{C}$-linear space $H_{n,s-1}^y$, 
provided that $x\in(\mathbb{C}^\ast)^s$ is generic. 
\end{lem}
\proof 
Since the dimension of $H_{n,s-1}^y$ is ${n+s-1\choose n}$, it suffices to show 
$\{F_\mu(x;y)\,;\, \mu\in Z_{s,n}\}$ is linearly independent. 
It is confirmed from the fact 
$$
\det F
=\det\big(F_\mu(x;\eta_{\nu}(x))\big)_{\mu,\nu\in Z}
=\prod_{\mu\in Z_{s,n}}F_\mu(x;\eta_{\mu}(x))\ne 0,
$$
which is a consequence of Lemma \ref{lem:triangular-F}.\qed

\par\medskip

\noindent
{\bf Proof of Theorem \ref{thm:E}.} 
We fix a generic point $x$ in $(\mathbb{C}^\ast)^s$.
Since $\Psi(z;y)$ 
as a function of $y$ is in $H_{n,s-1}^y$ by definition, using Lemma \ref{lem:basis-F},  
it is expressed as a linear combination of $F_\mu(x;y)$ ($\mu\in Z_{s,n}$), i.e. 
\begin{equation}
\label{eq:Psi=fF}
\Psi(z;y)=\prod_{i=1}^n\prod_{j=1}^{s-1}e(z_i;y_j)
=\sum_{\mu\in Z_{s,n}}f_\mu(z)F_\mu(x;y), 
\end{equation}
where $f_\lambda(z)$ are coefficients independent of $y$.  
Substituting $y=\eta_\nu(x)$ in this formula, we have 
$$
\Psi(z;\eta_\nu(x))
=\sum_{\mu\in Z_{s,n}}f_\mu(z)F_\mu(x;\eta_{\nu}(x))
\qquad (\nu \in Z_{s,n}).
$$
In what follows, 
we denote by $G=(G_{\mu,\nu}(x))_{\mu,\nu\in Z}$ 
the inverse matrix of 
$F=(F_{\mu}(x;\eta_\nu(x)))_{\mu,\nu\in Z}$.  
Then we obtain 
\begin{equation}
\label{eq:f in H}
f_\lambda(z)=
\sum_{\nu\in Z_{s,n}}
\Psi(z;\eta_\nu(x))G_{\nu\lambda}(x)
\qquad (\lambda \in Z_{s,n}), 
\end{equation}
which implies $f_\lambda(z)\in H_{s-1,n}^z$ 
since 
$\Psi(z;\eta_\nu(x))\in H_{s-1,n}^z$. 
Setting $z=x_\mu$ for each $\mu\in Z_{s,n}$
as in \eqref{eq:x_mu}, we obtain 
\begin{equation}
\label{eq:f=delta}
f_\lambda(x_\mu)
=\sum_{\nu\in Z_{s,n}}
\Psi(x_\mu,\eta_\nu(x))G_{\nu\lambda}(x)
=\sum_{\nu\in Z_{s,n}}F_\mu(x;\eta_\nu(x))G_{\nu\lambda}(x)
=\delta_{\lambda\mu}
\quad(\lambda,\mu\in Z_{s,n}).  
\end{equation}
These functions $f_\lambda(z)\in H^{z}_{s-1,n}$ ($\lambda\in Z_{s,n}$) 
are exactly what we wanted to construct.  
\qed

\par\medskip
Denoting the functions $f_\lambda(z)$ by 
$E_\lambda(x;z)$, we call them the  
{\em elliptic Lagrange 
interpolation functions} of type $BC_n$. 
We restate (\ref{eq:Psi=fF}) as corollary, 
which indicates the duality between 
$E_\lambda(x;z)$ and 
$F_\lambda(x;y)$.  
\begin{cor}[Duality]
\begin{equation}
\label{eq:Psi=EF}
\Psi(z;y)=\prod_{i=1}^n\prod_{j=1}^{s-1}e(z_i,y_j)
=\sum_{\lambda\in Z_{s,n}}E_\lambda(x;z)F_\lambda(x;y).
\end{equation}
\end{cor}
We will use this formula 
again in the succeeding sections. 

\par\medskip

From (\ref{eq:f in H}) we have 
\begin{equation}
\label{eq:E=GPsi}
E_\lambda(x;z)=
\sum_{\mu\in Z_{s,n}\atop \mu\preceq\lambda}
\Psi(z;\eta_\mu(x))G_{\mu\lambda}(x),
\end{equation}
which leads us to another explicit expression of 
$E_\lambda(x;z)$ different from (\ref{eq:E-explicit}) in Theorem \ref{thm:E-explicit}.
Recall that, for an invertible upper triangular matrix $A=(a_{ij})_{i,j=1}^{N}$  with $a_{ij}=0$ ($i>j$), the entries of 
its inverse 
are given by 
\begin{equation}\label{eq:Neumann}
\begin{split}
(A^{-1})_{ij}=\sum_{r=0}^{N-1}(-1)^r
\sum_{i=k_0<k_1<\cdots<k_r=j}
\frac{a_{k_0k_1}a_{k_1k_2}\cdots a_{k_{r-1}k_r}}
{a_{k_0k_0}a_{k_1k_1}a_{k_2k_2}\cdots a_{k_rk_r}}
\quad(i\le j). 
\end{split}
\end{equation}
In fact, setting  
$D=\mbox{diag}(a_{11},\ldots,a_{NN})$ 
and $B=A-D$, 
from $A=D(I+D^{-1}B)$ 
we obtain
\begin{equation*}
\begin{split}
A^{-1}&=(I+D^{-1}B)^{-1}D^{-1}
=\sum_{r=0}^{N-1}(-1)^r (D^{-1}B)^r D^{-1}\\
&=\sum_{r=0}^{N-1}(-1)^r 
D^{-1}BD^{-1}\cdots D^{-1}BD^{-1}, 
\end{split}
\end{equation*}
whose $(i,j)$-components are given by \eqref{eq:Neumann}.  
\begin{cor}
\label{cor:E-explicit-0}
The elliptic Lagrange interpolation functions 
$E_\lambda(x;z)$ are expressed explicitly as 
\begin{equation*}
E_\lambda(x;z)=\sum_{\mu\in Z_{s,n}\atop \mu\preceq\lambda}
\!\!\left(
\sum_{r\ge 0}(-1)^r
\!\!\!\!\!\!
\sum_{\mu=\nu^{(0)}\prec\cdots\prec\nu^{(r)}= \lambda}
\!
\frac
{\prod_{k=1}^r \prod_{i=1}^s\prod_{j=1}^{s-1}
e(x_i;x_j t^{\nu_j^{(k)}})_{\nu_i^{(k-1)}}}
{\prod_{k=0}^r \prod_{i=1}^s\prod_{j=1}^{s-1}e(x_i;x_jt^{\nu_j^{(k)}})_{\nu_i^{(k)}}}
\right)\!
\prod_{i=1}^n\prod_{j=1}^{s-1}e(z_i;x_jt^{\mu_j}).
\end{equation*}
\end{cor}
\proof 
In (\ref{eq:E=GPsi}), 
$G=\big(G_{\mu\nu}(x)\big)_{\mu,\nu\in Z}$ 
is the inverse matrix of the upper triangular matrix 
$F=\big(F_\mu(x;\eta_{\nu}(x))\big)_{\mu,\nu\in Z}$.  
By \eqref{eq:Neumann} we have the expression 
\begin{equation}
\begin{split}
G_{\mu\lambda}(x)
&=
\sum_{r\ge 0}(-1)^r
\sum_{\mu=\nu^{(0)}\prec\cdots\prec\nu^{(r)}= \lambda}
\frac{\prod_{k=1}^r F_{\nu^{(k-1)}}(x;\eta_{\nu^{(k)}})}{\prod_{k=0}^r F_{\nu^{(k)}}(x;\eta_{\nu^{(k)}})}
\\
&=\sum_{r\ge 0}(-1)^r
\sum_{\mu=\nu^{(0)}\prec\cdots\prec\nu^{(r)}= \lambda}
\frac{\prod_{k=1}^r \prod_{i=1}^s\prod_{j=1}^{s-1}e(x_i;x_jt^{\nu_j^{(k)}})_{\nu_i^{(k-1)}}}
{\prod_{k=0}^r \prod_{i=1}^s\prod_{j=1}^{s-1}e(x_i;x_jt^{\nu_j^{(k)}})_{\nu_i^{(k)}}}.
\label{eq:G-entries}
\end{split}
\end{equation}
By definition we also have 
\begin{equation}
\label{eq:Psi-eta}
\Psi(z;\eta_\mu(x))=
\prod_{i=1}^n\prod_{j=1}^{s-1}e(z_i;x_jt^{\mu_j}).
\end{equation}
Putting (\ref{eq:G-entries}) and (\ref{eq:Psi-eta}) on (\ref{eq:E=GPsi}) we have the expression in Corollary. \qed
\par\medskip 
\noindent
{\bf Remark.} The expression of $E_\lambda(x;z)$ in Corollary \ref{cor:E-explicit-0} 
is much more complex than (\ref{eq:E-explicit}) in Theorem \ref{thm:E-explicit}. 
Actually, it often becomes a huge sum 
even in the case where  
$E_\lambda(x;z)$ can be written simply 
in the form of product like (\ref{eq:ne_i}). 

%%%%%%%%%%%%%%%%%%%%%%%%%%%%%%%%%%
%%%%%%%%%%%%%%%%%%%%%%%%%%%%%%%%%%
\section{\boldmath 
Explicit expression %of the interpolation functions
for $E_\lambda(x;z)$}
%%%%%%%%%%%%%%%%%%%%%%%%%%%%%%%%%%
\label{Section:3}
%%%%%%%%%%%%%%%%%%%%%%%%%%%%%%%%%%
%%%%%%%%%%%%%%%%%%%%%%%%%%%%%%%%%%
In this section we give a proof of Theorem \ref{thm:E-explicit}. 
The main part of the proof is due to the repeated use of the recurrence relation for our interpolation functions. 
We first mention 
the explicit form of the interpolation 
functions of the case $n=1$ as the initial step of the recursive process.   
\begin{lem} 
For $x=(x_1,x_2,\ldots,x_s)\in (\mathbb{C}^*)^s$ and $z\in \mathbb{C}^*$
the interpolation functions 
of the case $n=1$ are expressed explicitly as 
\begin{equation}
\label{eq:E-explicit,n=1(2)}
E_{\epsilon_i}(x;z)=\prod_{1\le j\le s\atop j\ne i}
\frac{e(z;x_j)}{e(x_i;x_j)}
\end{equation}
\vskip -15pt
\noindent
for $i=1,2,\ldots,s$, where $\epsilon_i=(0,\ldots,0,\overset{\text{\tiny$i$}}{\stackrel{\text{\tiny$\smile$}}{1}},0,\ldots,0)
\in Z_{s,1}=\{\epsilon_1,\epsilon_2,\ldots,\epsilon_s\}$. 
\end{lem}
\proof Since $x_{\epsilon_i}=x_i$, it is immediately verified 
that the functions of the right-hand side of (\ref{eq:E-explicit,n=1(2)}) 
satisfy the conditions $E_{\epsilon_i}(x;z)\in H_{s-1,1}^{z}$ 
and $E_{\epsilon_i}(x;x_{\epsilon_j})=\delta_{ij}$. Such functions are determined uniquely by Theorem \ref{thm:E}. \qed 
\vskip 5pt 
Next we state the recursion formula for the interpolation 
functions.
\begin{lem}[Recursion formula] 
\label{lem:split-E}
Suppose that $n=m+l$. 
For $z=(z_1,\ldots,z_n)\in (\mathbb{C}^*)^n$ written as 
$z=(z',z'')$ where 
$z'=(z_1,\ldots,z_m)\in (\mathbb{C}^*)^m$ and $z''=(z_{m+1},\ldots,z_n)\in (\mathbb{C}^*)^l$, 
the function $E_\lambda(x;z)$  $(\lambda\in Z_{s,n})$ 
is expressed as 
\begin{equation}
\label{eq:split-E}
E_\lambda(x;z)=\sum_{{\mu\in Z_{s,m},\nu\in Z_{s,l}}\atop \mu+\nu=\lambda}
E_\mu(x;z')E_\nu(xt^\mu;z''),
\end{equation}
where
$xt^\mu=(x_1t^{\mu_1},x_2t^{\mu_2},\ldots,x_s t^{\mu_s})$ for $x=(x_1,x_2,\ldots,x_s)
\in (\mathbb{C}^*)^s$.
\end{lem}
\proof 
From the definition (\ref{eq:def-F}) of $F_\mu(x;y)$, it follows that 
$$
F_\mu(x;y)F_\nu(xt^\mu;y)=F_{\mu+\nu}(x;y).
$$
From this fact and (\ref{eq:Psi=EF}) we have 
\begin{eqnarray}
&&\hspace{-37pt}
\Psi(z;y)=\Psi(z';y)\Psi(z'';y)\nonumber\\[5pt]
&&=
\sum_{\mu\in Z_{s,m}}E_\mu(x;z')F_\mu(x;y)
\Psi(z'';y)
\nonumber\\
&&=
\sum_{\mu\in Z_{s,m}}E_\mu(x;z')F_\mu(x;y)
\Big(\sum_{\nu\in Z_{s,l}}E_\nu(xt^\mu;z'')F_\nu(xt^\mu;y)\Big)
\nonumber\\
&&=\sum_{\mu\in Z_{s,m}}\sum_{\nu\in Z_{s,l}}E_\mu(x;z')E_\nu(xt^\mu;z'')F_\mu(x;y)F_\nu(xt^\mu;y)
\nonumber\\
&&=\sum_{\mu\in Z_{s,m}}\sum_{\nu\in Z_{s,l}}E_\mu(x;z')E_\nu(xt^\mu;z'')F_{\mu+\nu}(x;y)
\nonumber\\
&&=\sum_{\lambda\in Z_{s,n}}\bigg[
\sum_{{\mu\in Z_{s,m},\nu\in Z_{s,l}}\atop \mu+\nu=\lambda}
E_\mu(x;z')E_\nu(xt^\mu;z'')
\bigg]
F_{\lambda}(x;y).
\label{eq:Psi=EEF}
\end{eqnarray}
Comparing (\ref{eq:Psi=EEF}) with (\ref{eq:Psi=EF}), we obtain the expression (\ref{eq:split-E}) in Lemma \ref{lem:split-E}.\qed
\par\medskip
\noindent
{\bf Remark.} 
The special cases $m=1$ or $l=1$ of the recursion formula in Lemma \ref{lem:split-E} indicate that 
\begin{equation}
\label{eq:split-E n=1(1)}
E_\lambda(x;z)=\sum_{i=1}^s
E_{\epsilon_{i}}(x;z_1)E_{\lambda-\epsilon_{i}}(xt^{\epsilon_{i}};z_2,\ldots,z_n)
\end{equation}
or 
\begin{equation}
\label{eq:split-E n=1(2)}
E_\lambda(x;z)=\sum_{i=1}^s
E_{\lambda-\epsilon_{i}}(x;z_1,\ldots,z_{n-1})E_{\epsilon_{i}}(xt^{\lambda-\epsilon_{i}};z_n),
\end{equation}
where we regard $E_{\lambda-\epsilon_{i}}(x;z_1,\ldots,z_{n-1})=0$ if $\lambda-\epsilon_{i}\not\in Z_{s,n-1}$.\\

By the repeated use of (\ref{eq:split-E n=1(1)}) or (\ref{eq:split-E n=1(2)}) we have the following expression.
\begin{cor} 
\label{cor:E-explicit}
$$
E_\lambda(x;z)=\sum_{{(i_1,\ldots,i_n)\in \{1,\ldots,s\}^n}\atop\epsilon_{i_1}+\cdots+\epsilon_{i_n}=\lambda}
E_{\epsilon_{i_1}}(x;z_1)
E_{\epsilon_{i_2}}(xt^{\epsilon_{i_1}};z_2)
E_{\epsilon_{i_3}}(xt^{\epsilon_{i_1}+\epsilon_{i_2}};z_3)
\cdots
E_{\epsilon_{i_n}}(xt^{\epsilon_{i_1}+\cdots+\epsilon_{i_{n-1}}};z_n).
$$
\end{cor}
\par\medskip
Rewriting 
Corollary \ref{cor:E-explicit}, 
we obtain the explicit formula (\ref{eq:E-explicit}) for 
$E_\lambda(x;z)$ as presented 
in Theorem \ref{thm:E-explicit}. 
\begin{thm} 
The interpolation functions 
$E_\lambda(x;z)$ are expressed explicitly as 
\begin{equation}
\label{eq:E-explicit-2}
E_\lambda(x;z)=\sum_{{(i_1,\ldots,i_n)\in \{1,\ldots,s\}^n}\atop\epsilon_{i_1}+\cdots+\epsilon_{i_n}=\lambda}
\prod_{k=1}^n
\prod_{1\le j\le s\atop j\ne i_k}\frac
{e(z_k\,;x_jt^{\lambda_j^{(k-1)}})}
{e(x_{i_k}t^{\lambda_{i_k}^{(k-1)}};x_jt^{\lambda_j^{(k-1)}})},
\end{equation}
where $\lambda_i^{(k)}=|\{l\in\{1,\ldots,k\}\,;\,i_l=i\}|$. 
Equivalently $E_\lambda(x;z)$ is also written as {\rm (\ref{eq:E-explicit})} in Theorem \ref{thm:E-explicit}, i.e., 
\begin{equation}
\label{eq:E-explicit-3}
E_\lambda(x;z)=\sum_{K_1\sqcup\cdots\sqcup K_s\atop =\{1,2,\ldots,n\}}
\prod_{i=1}^s\prod_{k\in K_i}
\prod_{1\le j\le s\atop j\ne i}
\frac
{e(z_k\,;x_jt^{\lambda_j^{(k-1)}})}
{e(x_it^{\lambda_i^{(k-1)}};x_jt^{\lambda_j^{(k-1)}})},
\end{equation}
where $\lambda_i^{(k)}=|K_i\cap\{1,2,\ldots,k\}|$ and the summation is taken over all index sets $K_i$ $(i=1,2,\ldots,s)$ satisfying 
$|K_i|=\lambda_i$ and $K_1\sqcup\cdots\sqcup K_s=\{1,2,\ldots,n\}$. 
\end{thm}

\proof
For $(i_1,\ldots,i_n)\in \{1,\ldots,s\}^n$ satisfying 
$\epsilon_{i_1}+\cdots+\epsilon_{i_n}=\lambda$, 
we set $\lambda^{(k)}=\epsilon_{i_1}+\cdots+\epsilon_{i_k}$
for $k=0,1,\ldots,n$. 
Then by definition $\lambda^{(k)}=(\lambda_1^{(k)},\ldots, \lambda_s^{(k)})\in Z_{s,k}$ is expressed by 
$$
\lambda_i^{(k)}=|\{l\in\{1,\ldots,k\}\,;\,i_l=i\}|\quad(i=1,\ldots,s).
$$
From Corollary \ref{cor:E-explicit}, we therefore obtain 
\begin{equation}
\label{eq:E-explicit-4}
E_\lambda(x;z)=\sum_{{(i_1,\ldots,i_n)\in \{1,\ldots,s\}^n}\atop\epsilon_{i_1}+\cdots+\epsilon_{i_n}=\lambda}
\prod_{k=1}^n E_{\epsilon_{i_k}}(xt^{\lambda^{(k-1)}};z_k), 
\end{equation}
which coincides with (\ref{eq:E-explicit-2}) using (\ref{eq:E-explicit,n=1(2)}).

Next we explain the latter part of the theorem.
Let $K_i$ be sets of indices specified by 
$K_i=\{l\in\{1,\ldots,n\}\,;\,i_l=i\}$, 
where $(i_1,\ldots,i_n)\in \{1,\ldots,s\}^n$ and 
$\epsilon_{i_1}+\cdots+\epsilon_{i_n}=\lambda$. 
Then $\lambda_i^{(k)}$ is written as 
$$
\lambda_i^{(k)}=|K_i\cap \{1,2,\ldots, k\}|.  
$$
In particular, we have 
$
\lambda_i=\lambda_i^{(n)}=|K_i|.
$
Thus $K_i$ $(i=1,2,\ldots,s)$ satisfy $K_1\sqcup\cdots\sqcup K_s=\{1,2,\ldots,n\}.$
Since $i_k=i$ if and only if $k\in K_i$, the expression (\ref{eq:E-explicit-4}) is rewritten as 
\begin{equation}
\label{eq:E-explicit-5}
E_\lambda(x;z)=\sum_{K_1\sqcup\cdots\sqcup K_s\atop =\{1,2,\ldots,n\}}
\prod_{i=1}^s\prod_{k\in K_i}E_{\epsilon_{i}}(xt^{\lambda^{(k-1)}};z_k)
\end{equation}
where the summation is taken over all index sets $K_i$ satisfying 
$|K_i|=\lambda_i$ and $K_1\sqcup\cdots\sqcup K_s=\{1,2,\ldots,n\}$. 
Therefore (\ref{eq:E-explicit-5}) coincides with (\ref{eq:E-explicit-3}) using (\ref{eq:E-explicit,n=1(2)}).\qed
\par\medskip
We remark that 
the interpolation functions of the special cases 
$\lambda=n\epsilon_i\in Z_{s,n}$ have simple 
factorized forms; this fact will be used 
in the succeeding section. 
\begin{cor} 
\label{cor:ne_i}
For $n\epsilon_i\in Z_{s,n}$ $(i=1,\ldots, s)$, 
one has 
\begin{equation}
\label{eq:ne_i}
E_{n\epsilon_i}(x;z)
=\prod_{1\le j\le s\atop j\ne i}
\frac{e(z_1;x_j)\cdots e(z_n;x_j)}{e(x_i;x_j)_n}. 
\end{equation}
\end{cor} 
\proof If we put $\lambda=n\epsilon_i$ in the formula of  Corollary \ref{cor:E-explicit}, 
then the right-hand side reduces to a single term 
with $(i_1,i_2,\ldots,i_n)=(i,i,\ldots,i)$. 
Therefore, using (\ref{eq:E-explicit,n=1(2)}) we obtain 
$$
E_{n\epsilon_i}(x;z)
=E_{\epsilon_{i}}(x;z_1)
E_{\epsilon_{i}}(xt^{\epsilon_{i}};z_2)
E_{\epsilon_{i}}(xt^{2\epsilon_{i}};z_3)
\cdots
E_{\epsilon_{i}}(xt^{(n-1)\epsilon_{i}};z_n)
$$
$$
=\prod_{1\le j\le s\atop j\ne i}\frac{e(z_1;x_j)}{e(x_i;x_j)}\frac{e(z_2;x_j)}{e(x_it;x_j)}\frac{e(z_3;x_j)}{e(x_it^2;x_j)}\cdots\frac{e(z_n;x_j)}{e(x_it^{n-1};x_j)},
$$
which coincides with (\ref{eq:ne_i}).\qed

%%%%%%%%%%%%%%%%%%%%%%%%%%%%%%%%%%
%%%%%%%%%%%%%%%%%%%%%%%%%%%%%%%%%%
\section{Transition coefficients for the interpolation functions}
%%%%%%%%%%%%%%%%%%%%%%%%%%%%%%%%%%
\label{Section:4}
%%%%%%%%%%%%%%%%%%%%%%%%%%%%%%%%%%
%%%%%%%%%%%%%%%%%%%%%%%%%%%%%%%%%%
In this section we 
discuss the transition coefficients between two sets of  interpolation functions 
with different parameters. 

For generic $x,y\in (\mathbb{C}^*)^s$, the interpolation functions $E_\mu(x;z)\in H_{s-1,n}^z$ as functions of 
$z\in (\mathbb{C}^*)^n$ are expanded in terms of $E_\nu(y;z)$ ($\nu\in Z_{s,n}$), i.e., 
\begin{equation}
\label{eq:connection-E}
E_\mu(x;z)=\sum_{\nu\in Z_{s,n}} C_{\mu\nu}(x;y)E_\nu(y;z),
\end{equation}
where the coefficients $C_{\mu\nu}(x;y)$ are independent of 
$z$.  
From the property (\ref{eq:E-delta}) of the 
interpolation functions, 
we immediately see that $C_{\mu\nu}(x;y)$ is expressed by the special value of $E_\mu(x;z)$ as
$$
C_{\mu\nu}(x;y)=E_\mu (x;y_{\nu})\quad(\mu,\nu\in Z_{s,n}). 
$$ 
For $x,y\in (\mathbb{C}^*)^s$, 
we denote the {\em transition matrix} from 
$(E_\lambda(x;z))_{\lambda\in Z_{s,n}}$ 
to 
$(E_\lambda(y;z))_{\lambda\in Z_{s,n}}$
by 
$$
E(x;y)=\Big(E_\mu (x;y_{\nu})\Big)_{\!\mu,\nu\in Z_{s,n}},
$$
where the rows and the columns 
are arranged in the total order $\prec$ of $Z_{s,n}$. 
By definition, for generic 
$x,y,w\in (\mathbb{C}^*)^s$ we have 
\begin{equation}
\label{eq:E=EE}
E(x;y)=E(x;w)E(w;y),
\end{equation}
in particular
\begin{equation}
E(x;x)=I
\quad
\mbox{and} 
\quad
E(y;x)=
E(x;y)^{-1}.
\end{equation} 
\begin{thm}
\label{thm:detE}
For generic $x,y\in (\mathbb{C}^*)^s$ the determinant of 
the transition matrix 
$E(x;y)$ is given explicitly by 
\begin{equation}
\label{eq:detE-1}
\det E(x;y)
=\prod_{k=1}^n\bigg[\prod_{r=0}^{n-k}\prod_{1\le i<j\le s}
\frac{e(y_it^r;y_jt^{(n-k)-r})}{e(x_it^r;x_jt^{(n-k)-r})}
\bigg]^{s+k-3\choose k-1},
\end{equation}
or equivalently by 
\begin{equation}
\label{eq:detE-2}
\det E(x;y)
=
\prod_{k=1}^n\bigg[
\prod_{r=0}^{n-k}\prod_{1\le i<j\le s}
\frac
{x_i\theta(t^{2r-(n-k)}y_iy_j^{-1})\theta(t^{n-k}y_iy_j)}
{y_i\theta(t^{2r-(n-k)}x_ix_j^{-1})\theta(t^{n-k}x_ix_j)}
\bigg]^{{s+k-3\choose k-1}}.
\end{equation}
\end{thm}

The goal of this section is to prove the above theorem. 
For this purpose we first investigate a special case. 
\begin{lem}
\label{lem:E-triangular}
For $x,y\in (\mathbb{C}^*)^s$ suppose that $y_i=x_i$ $(i=1,2,\ldots,s-1)$, i.e. 
$y=(x_1,\ldots,x_{s-1},y_s)$. 
For $\alpha,\beta\in Z_{s,n}$ if there exists $i\in\{1,2,\ldots,s-1\}$ such that $\alpha_i<\beta_i$, then 
$
E_\alpha(x;y_\beta)=0
$. In particular, $E(x,y)$ is a lower 
triangular matrix with the diagonal entries 
$$
E_\alpha(x;y_\alpha)
=\prod_{i=1}^{s-1}\frac{e(y_s;x_it^{\alpha_i})_{\alpha_s}}{e(x_s;x_it^{\alpha_i})_{\alpha_s}}.
$$
Moreover the determinant of $E(x,y)$ of the case $y=(x_1,\ldots,x_{s-1},y_s)$ is expressed as 
\begin{equation}
\label{eq:detE-3}
\det E(x,y)=
\prod_{\alpha\in Z_{s,n}}\prod_{i=1}^{s-1}
\frac{e(y_s;x_it^{\alpha_i})_{\alpha_s}}{e(x_s;x_it^{\alpha_i})_{\alpha_s}}
=\prod_{k=1}^n\bigg[\prod_{r=0}^{n-k}\prod_{i=1}^{s-1}
\frac{e(x_it^r;y_st^{(n-k)-r})}{e(x_it^r;x_st^{(n-k)-r})}
\bigg]^{s+k-3\choose k-1}.
\end{equation}
\end{lem}
\proof
If $\beta_s=0$ for $\beta\in Z_{s,n}$, then 
$C_{\alpha\beta}(x,y)=E_\alpha(x;y_\beta)=\delta_{\alpha\beta}$ 
by the definition (\ref{eq:connection-E}).
If $\beta_s\ne 0$ for $\beta\in Z_{s,n}$, 
we apply Lemma \ref{lem:split-E} 
with $m=\beta_1+\cdots+\beta_{s-1}$, $l=\beta_s$ 
to obtain 
\begin{equation}
\label{eq:C_ab-1}
\begin{split}
C_{\alpha\beta}(x;y)&=E_\alpha(x;y_\beta)
=\sum_{{\mu\in Z_{s,m},\nu\in Z_{s,l}}\atop \mu+\nu=\alpha}
E_\mu(x;y_{\beta'})E_\nu(xt^\mu;y_{\beta''})\\
&=\sum_{{\mu\in Z_{s,m},\nu\in Z_{s,l}}\atop \mu+\nu=\alpha}
E_\mu(x;x_{\beta'})E_\nu(xt^\mu;y_{\beta''}),
\end{split}
\end{equation}
where
$\beta'=(\beta_1,\ldots,\beta_{s-1},0)\in Z_{s,m}$, 
and $\beta''=(0,\ldots,0,\beta_s)\in Z_{s,l}$ 
so that $y_{\beta'}=x_{\beta'}$.  
Note that 
$$
x_{\beta'}=(
\underbrace{x_1,x_1t,\ldots,x_1t^{\beta_1-1}\phantom{\Big|}\!\!}_{\beta_1}
,\underbrace{x_2,x_2t,\ldots,x_2t^{\beta_2-1}\phantom{\Big|}\!\!}_{\beta_2}
,\ldots
,\underbrace{x_{s-1},x_{s-1}t,\ldots,x_{s-1}t^{\beta_{s-1}-1}\phantom{\Big|}\!\!}_{\beta_{s-1}})
\in (\mathbb{C}^*)^m,
$$
$$
y_{\beta''}
=(y_s,y_st,\ldots,y_st^{\beta_s-1})\in (\mathbb{C}^*)^l.
$$
By the property of the interpolation functions, 
we have $E_{\mu}(x;x_{\beta'})=\delta_{\mu\beta'}$
for $\mu\in Z_{s,m}$. 
From (\ref{eq:C_ab-1}) $C_{\alpha\beta}(x;y)$ is written as 
\begin{equation}
\label{eq:C_ab-2}
C_{\alpha\beta}(x;y)=\sum_{{\mu\in Z_{s,m},\nu\in Z_{s,l}}\atop \mu+\nu=\alpha}
\delta_{\mu\beta'}E_\nu(xt^\mu;y_{\beta''})
=E_{\alpha-\beta'}(xt^{\beta'};y_{\beta''}),
\end{equation}
where $\alpha-\beta'=(\alpha_1-\beta_1,\ldots,\alpha_{s-1}-\beta_{s-1}, \alpha_s)\in Z_{s,l}$ 
and 
$xt^{\beta'}=(x_1t^{\beta_1},x_2t^{\beta_2},\ldots,x_{s-1} t^{\beta_{s-1}},x_s)
\in (\mathbb{C}^*)^s$. 
If there exists $i\in \{1,2,\ldots,s-1\}$ such that $\alpha_i<\beta_i$, i.e., $\alpha-\beta'\not\in Z_{s,l}$, 
then $C_{\alpha\beta}(x;y)=0$. 
In particular, if $\alpha\prec\beta$, then $C_{\alpha\beta}(x;y)=0$, 
which indicates the matrix $E(x,y)$ of the case $y=(x_1,\cdots,x_{s-1},y_s)$ is lower triangular. 
On the other hand, if 
$\alpha_1\ge\beta_1$, $\alpha_2\ge\beta_2, \ldots, \alpha_{s-1}\ge\beta_{s-1}$
(and $\alpha_s=n-(\alpha_1+\cdots+\alpha_{s-1})\le n-(\beta_1+\cdots+\beta_{s-1})=\beta_s$), 
$
C_{\alpha\beta}(x;y)
$ is written as (\ref{eq:C_ab-2}).
In particular, if $\alpha=\beta$, then 
\begin{equation}
\label{eq:C_bb}
C_{\beta\beta}(x;y)=E_{\beta''}(xt^{\beta'};y_{\beta''}). 
\end{equation}
Since $\beta''=\beta_s \epsilon_s$, 
from (\ref{eq:ne_i}) in Corollary \ref{cor:ne_i} we have 
\begin{equation}
\label{eq:E_beta_s}
E_{\beta''}(xt^{\beta'};z_{m+1},\ldots,z_n)
=\prod_{i=1}^{s-1}\frac{e(z_{m+1};x_it^{\beta_i})\cdots e(z_n;x_it^{\beta_i})}
{e(x_s;x_it^{\beta_i})_{\beta_s}}.
\end{equation}
From (\ref{eq:C_bb}) and (\ref{eq:E_beta_s}) we therefore obtain 
$$
C_{\beta\beta}(x;y)=
E_\beta(x;y_\beta)
=\prod_{i=1}^{s-1}\frac{e(y_s;x_it^{\beta_i})_{\beta_s}}
{e(x_s;x_it^{\beta_i})_{\beta_s}}.
$$

Lastly we derive \eqref{eq:detE-3}.  
Since the matrix $E(x;y)$ is lower triangular, its determinant 
is calculated as
\begin{equation}
\begin{split}
&\det\Big(E_\alpha(x;y_\beta)\Big)_{\alpha,\beta\in Z_{s,n}}
=
\prod_{\alpha\in Z_{s,n}}\prod_{i=1}^{s-1}
\frac{e(y_s;x_it^{\alpha_i})_{\alpha_s}}{e(x_s;x_it^{\alpha_i})_{\alpha_s}}
\\
&=
\prod_{i=1}^{s-1}
\prod_{\alpha\in Z_{s,n}}
\prod_{l=1}^{\alpha_s}
\frac{e(y_st^{l-1};x_it^{\alpha_i})}
{e(x_st^{l-1};x_it^{\alpha_i})}
=
\prod_{i=1}^{s-1}
\prod_{r=0}^{n-1}
\prod_{l=1}^{n-r}
\prod_{\alpha\in Z_{s,n}\atop
\,\alpha_i=r,\,\alpha_s\ge l}
\frac{e(y_st^{l-1};x_it^r)}
{e(x_st^{l-1};x_it^r)}. 
\end{split}
\end{equation}
Here we count the number of $\alpha\in Z_{s,n}$ 
such that $\alpha_i=r$ and $\alpha_s\ge l$.  
Note that for $0\le r+k\le n$, 
$|\{\alpha\in Z_{s,n}\,;\,\alpha_i=r, \alpha_s=k\}|
=|Z_{s-2,n-r-k}|
={\textstyle \binom{n-r-k+s-3}{s-3}}. 
$
Hence for $r+l\le n$, 
we have 
$|\{\alpha\in Z_{s,n}\,;\,\alpha_i=r, \alpha_s\ge l\}|
=\sum_{k=l}^{n-r}{\textstyle \binom{n-r-k+s-3}{s-3}}
=\sum_{p=0}^{n-r-l}{\textstyle\binom{p+s-3}{s-3}}
=
{\textstyle \binom{n-r-l+s-2}{s-2}}.
$
Finally we obtain
\begin{equation}
\begin{split}
&\det\Big(E_\alpha(x;y_\beta)\Big)_{\alpha,\beta\in Z_{s,n}}
=
\prod_{i=1}^{s-1}
\prod_{r=0}^{n-1}
\prod_{l=1}^{n-r}
\left(
\frac{e(y_st^{l-1};x_it^r)}
{e(x_st^{l-1};x_it^r)}
\right)^{\binom{n-r-l+s-2}{s-2}}
\\
&=
\prod_{i=1}^{s-1}
\prod_{r=0}^{n-1}
\prod_{k=1}^{n-r}
\left(
\frac{e(y_st^{n-r-k};x_it^r)}
{e(x_st^{n-r-k};x_it^r)}
\right)^{\binom{k+s-3}{s-2}}
%\\
%&
=
\prod_{i=1}^{s-1}
\prod_{k=1}^{n}
\prod_{r=0}^{n-k}
\left(
\frac{e(y_st^{n-r-k};x_it^r)}
{e(x_st^{n-r-k};x_it^r)}
\right)^{\binom{k+s-3}{k-1}},
\end{split}
\end{equation}
which coincides with (\ref{eq:detE-3}). The proof is now complete. \qed

\par\medskip
We now prove Theorem \ref{thm:detE}. \\[10pt]
\noindent
{\bf Proof of Theorem \ref{thm:detE}.}
We set 
$$w^{(i)}=(x_1,\ldots,x_{i},y_{i+1},\ldots,y_s)\in (\mathbb{C}^*)^s$$
for $i=0,1,\ldots, s$, which satisty  
$w^{(s)}=x$ and $w^{(0)}=y$. 
Since $w^{(s-1)}=(x_1,\ldots,x_{s-1},y_s)$, Lemma \ref{lem:E-triangular} indicates that 
$$
\det E(w^{(s)};w^{(s-1)})=
\prod_{\alpha\in Z_{s,n}}\prod_{i=1}^{s-1}
\frac{e(y_s;x_it^{\alpha_i})_{\alpha_s}}{e(x_s;x_it^{\alpha_i})_{\alpha_s}}
=\prod_{k=1}^n\bigg[\prod_{r=0}^{n-k}\prod_{i=1}^{s-1}
\frac{e(x_it^r;y_st^{(n-k)-r})}{e(x_it^r;x_st^{(n-k)-r})}
\bigg]^{s+k-3\choose k-1}.
$$
In the same way as Lemma \ref{lem:E-triangular}, 
for $l=1,\ldots,s$ 
we have 
\begin{eqnarray}
&&
\hspace{-20pt}
\det E(w^{(l)};w^{(l-1)})=
\prod_{\alpha\in Z_{s,n}}\bigg[
\bigg(
\prod_{1\le i<l}
\frac{e(y_l;x_it^{\alpha_i})_{\alpha_l}}{e(x_l;x_it^{\alpha_i})_{\alpha_l}}
\bigg)
\bigg(
\prod_{l<j\le s}
\frac{e(y_l;y_jt^{\alpha_j})_{\alpha_l}}{e(x_l;y_jt^{\alpha_j})_{\alpha_l}}
\bigg)
\bigg]
\nonumber\\
&&=
\prod_{k=1}^n\prod_{r=0}^{n-k}
\bigg[
\bigg(
\prod_{1\le i<l}
\frac{e(x_it^r;y_lt^{(n-k)-r})}{e(x_it^r;x_lt^{(n-k)-r})}
\bigg)
\bigg(
\prod_{l<j\le s}
\frac{e(y_lt^r;y_jt^{(n-k)-r})}{e(x_lt^r;y_jt^{(n-k)-r})}
\bigg)
\bigg]^{s+k-3\choose k-1}
\label{eq:E(ww)}
\end{eqnarray}
exchanging the roles of indices $l$ and $s$. 
From the relation (\ref{eq:E=EE}) 
of transition matrices, we have the decomposition of $E(x;y)$ as 
\begin{equation}
\label{eq:E=EE...E}
E(x;y)=E(w^{(s)};w^{(0)})=E(w^{(s)};w^{(s-1)})E(w^{(s-1)};w^{(s-2)})\cdots E(w^{(1)};w^{(0)}).
\end{equation}
Applying (\ref{eq:E(ww)}) to (\ref{eq:E=EE...E}) we obtain 
$$\det E(x;y)=\prod_{l=1}^s\det E(w^{(l-1)};w^{(l)})
=\prod_{k=1}^n\bigg[\prod_{r=0}^{n-k}\prod_{1\le i<j\le s}
\frac{e(y_it^r;y_jt^{(n-k)-r})}{e(x_it^r;x_jt^{(n-k)-r})}
\bigg]^{s+k-3\choose k-1},
$$
which completes the proof of Theorem \ref{thm:detE}. \qed

\par\medskip 
\noindent
{\bf Remark.}  
In the decomposition \eqref{eq:E=EE...E}, 
each component 
$E(w^{(l)};w^{(l-1)})$ ($l=1,\ldots,s$) is 
lower triangular with respect to the partial ordering  
$\subseteq_{l}$ of $Z_{s,n}$ 
defined by 
$\alpha\subseteq_l\beta\ \ \Longleftrightarrow\ \ 
\alpha_i\le \beta_i\ \ (i\ne l). $

%%%%%%%%%%%%%%%%%%%%%%%%%%%%%%%%%%
%%%%%%%%%%%%%%%%%%%%%%%%%%%%%%%%%%
\section{\boldmath
Proofs of the main theorems for $BC_n$ Jackson integrals}
%%%%%%%%%%%%%%%%%%%%%%%%%%%%%%%%%%
\label{Section:5}
%%%%%%%%%%%%%%%%%%%%%%%%%%%%%%%%%%
%%%%%%%%%%%%%%%%%%%%%%%%%%%%%%%%%%
We conclude this paper by providing with proofs of Theorems \ref{thm:Wronski} and \ref{thm:connection2}, 
on the basis of properties of the $BC_n$ elliptic 
Langrange interpolation functions $E_\lambda(x;z)$ 
as we established in the previous sections.

\par\medskip
\noindent
{\bf Proof of Theorem \ref{thm:connection2}.}
If  $\varphi(z)$ is a $W_n$-invariant 
holomorphic function on 
$(\mathbb{C}^\ast)^n$, 
then $\la \varphi,z\ra$ is a 
meromorphic function on 
$(\mathbb{C}^\ast)^n$ and $q$-periodic 
with respect to each variable $z_i$ ($i=1,\ldots,s$). 
It is known by \cite[Definition 3.8]{AI2009} that 
the regularization $\lla\varphi,z\rra=\la \varphi,z\ra/\Theta(z)$ 
is a $W_n$-invariant holomorphic function on 
$(\mathbb{C}^\ast)^n$, and belongs to 
$H_{s-1,n}^z$ due to the quasi-periodicity of $1/\Theta(z)$.
This implies that 
$\lla\varphi, z\rra$ is expressed as a linear combination of our 
interpolation functions $E_\mu(x;z)$ ($\mu\in Z_{s,n})$, i.e.,
$$
\lla\varphi, z\rra=\sum_{\mu\in Z}d_\mu E_\mu(x;z).
$$
From (\ref{eq:E-delta}), we obtain 
$d_\nu=\sum_{\mu\in Z}d_\mu E_\mu(x;x_\nu)=\lla\varphi, x_\nu\rra$, 
which completes the proof of Theorem \ref{thm:connection2}. \qed
\vspace{10pt}
\noindent
{\bf Proof of Theorem \ref{thm:Wronski}.}
%\proof 
The spacial case $x=a=(a_1,\ldots,a_s)$ of Theorem \ref{thm:Wronski} was proved 
in \cite[Theorem 1.3]{AI2009}, i.e.,   
\begin{eqnarray}
\label{eq:Wronski2}
&&
\hspace{-16pt}
\det\Big(\lla\chi_\lambda,a_\mu\rra\Big)_{\!\!\lambda\in B\atop \!\!\mu\in Z}
=
\prod_{k=1}^n\bigg[\bigg(\!
(1-q)\frac{(q)_\infty(qt^{-(n-k+1)})_\infty}{(qt^{-1})_\infty}\!\bigg)^{\!\! s}\,
\frac{\prod_{1\le i<j\le 2s+2}(qt^{-(n-k)}a_i^{-1}a_j^{-1})_\infty}{(qt^{-(n+k-2)}a_1^{-1}a_2^{-1}\cdots a_{2s+2}^{-1})_\infty}
\bigg]^{{s+k-2\choose k-1}}
\nonumber\\
&&
\hspace{90pt}
\times\prod_{k=1}^n\bigg[
\prod_{r=0}^{n-k}\prod_{1\le i<j\le s}
\frac{\theta(t^{2r-(n-k)}a_ia_j^{-1})\theta(t^{n-k}a_ia_j)}{t^ra_i}
\bigg]^{{s+k-3\choose k-1}}.
\end{eqnarray}
On the other hand, 
if we put $x=a=(a_1,\ldots,a_s)$ on (\ref{eq:connection2}) in Theorem \ref{thm:connection2}, 
then we have 
$$
\lla \varphi, z\rra=\sum_{\mu\in Z}\lla\varphi, a_\mu\rra E_\mu(a;z). 
$$ 
In particular, 
setting $\varphi (z)=\chi_\lambda(z)$ ($\lambda\in B_{s,n}$) and $z=x_\nu$ ($\nu\in Z_{s,n}$) we obtain
$$
\Big(\lla\chi_\lambda,x_\nu\rra\Big)_{\lambda\in B\atop \nu\in Z}
=
\Big(\lla\chi_\lambda,a_\mu\rra\Big)_{\lambda\in B\atop \mu\in Z}
\Big(E_\mu(a;x_\nu)\Big)_{\mu\in Z\atop \nu\in Z},
$$
so that 
$$\displaystyle
\det\Big(\lla\chi_\lambda,x_\nu\rra\Big)_{\lambda\in B\atop \nu\in Z}
=
\det\Big(\lla\chi_\lambda,a_\mu\rra\Big)_{\lambda\in B\atop \mu\in Z}
\det E(a;x).  
$$
Combining \eqref{eq:Wronski2} and 
\eqref{eq:detE-2}, 
we obtain the determinant formula of \eqref{eq:Wronski}.  
This completes the proof of Theorem \ref{thm:Wronski}.\qed

%%%%%%%%%%%%%%%%%%%%%%%%%%%%%%%%%%%%%%%%%
%%%%%%%%%%%%%%%%%%%%%%%%%%%%%%%%%%%%%%%%%
\section*{Acknowledgements} 
This work is supported by JSPS Kakenhi Grants (C)25400118 and 
(B)15H03626. 
%%%%%%%%%%%%%%%%%%%%%%%%%%%%%%%%%%%%%%%%%
%%%%%%%%%%%%%%%%%%%%%%%%%%%%%%%%%%%%%%%%%

\end{document}